\documentclass{amsart}
\usepackage{enumerate}
\usepackage{tabu}
\usepackage{longtable}
\usepackage{multirow}
\usepackage{amssymb} 
\usepackage{amsmath} 
\usepackage{amscd}
\usepackage{amsbsy}
\usepackage{comment}
\usepackage[matrix,arrow]{xy}
\usepackage{hyperref}

\DeclareMathOperator{\Frob}{Frob}

\DeclareMathOperator{\Gal}{Gal}

\DeclareMathOperator{\ord}{\upsilon}

\newcommand{\Q}{{\mathbb Q}}
\newcommand{\Z}{{\mathbb Z}}
\newcommand{\F}{{\mathbb F}}
\newcommand{\cA}{\mathcal{A}}

\newcommand{\mP}{\mathfrak{P}}

\def\mod#1{{\ifmmode\text{\rm\ (mod~$#1$)}
\else\discretionary{}{}{\hbox{ }}\rm(mod~$#1$)\fi}}

\begin {document}

\newtheorem{thm}{Theorem}
\newtheorem{lem}{Lemma}[section]
\newtheorem{prop}[lem]{Proposition}

\theoremstyle{definition}

\theoremstyle{remark}

\title[]{On powers\\ that are sums of consecutive like powers}


\author{Vandita Patel}
\address{Mathematics Institute, University of Warwick, Coventry CV4 7AL, United Kingdom}
\email{vandita.patel@warwick.ac.uk}

\author{Samir Siksek}
\address{Mathematics Institute, University of Warwick, Coventry CV4 7AL, United Kingdom}
\email{S.Siksek@warwick.ac.uk}
\thanks{
The first-named author is supported by an EPSRC studentship.
The second-named author
is supported by 
the EPSRC {\em LMF: L-Functions and Modular Forms} Programme Grant
EP/K034383/1.
}

\date{\today}

\keywords{Exponential equation, Bernoulli polynomial, Newton polygon}
\subjclass[2010]{Primary 11D61, Secondary 11B68}

\begin {abstract}
Let $k \ge 2$ be even, and let $r$ be a non-zero integer. 
We show that for almost all $d \ge 2$ (in the sense of natural density),
the equation
\[
x^k+(x+r)^k+\cdots+(x+(d-1)r)^k=y^n,\qquad x,~y,~n \in \Z, \qquad n \ge 2,
\]
has no solutions. 
\end {abstract}
\maketitle

\section{Introduction} \label{intro}

The problem of cubes that are sums of consecutive cubes goes
back to Euler
\cite[art.\ 249]{Euler} who noted the remarkable relation
$3^3+4^3+5^3=6^3$. Similar problems were considered by several 
mathematicians during the 19th and early 20th century
as surveyed in
Dickson's 
\emph{History of the Theory of Numbers}
\cite[pp. 582--588]{Dickson}. These questions are still of
interest today. 
For example,  
both Cassels \cite{Cassels}
and Uchiyama \cite{Uchiyama} determined
the squares that can be written as sums of three consecutive cubes.
Stroeker \cite{Stroeker}
determined all squares that are expressible as the sum of $2 \le d \le 50$
consecutive cubes, using a method based on linear forms in 
elliptic logarithms. More recently, Bennett, Patel and Siksek \cite{BPS2}
determined all perfect powers that are expressible as sums of
$2 \le d \le 50$ consecutive cubes, using linear forms in logarithms,
sieving and Frey curves. There has been some interest in powers
that are sums of $k$-th powers for other exponents $k$.
For example, the solutions to the equation 
\[
x^k+(x+1)^k+(x+2)^k=y^n, \qquad x,~y,~n \in \Z, \qquad n \ge 2,
\] 
have been determined by Zhongfeng Zhang for $k=2$, $3$, $4$
and by Bennett, Patel and Siksek \cite{BPS1} for $k=5$, $6$.

\medskip

In view of the above, it is natural to consider the equation
\begin{equation}\label{eqn:premain}
x^k+(x+1)^k+\cdots+(x+d-1)^k=y^n,\qquad x,~y,~n \in \Z, \qquad n \ge 2
\end{equation}
with $k$, $d \ge 2$.
This was studied by Zhang and Bai \cite{ZB} for $k=2$. They show that if 
$q$ is a prime $\equiv \pm 5 \pmod{12}$ and $\ord_q(d)=1$
then equation \eqref{eqn:premain} has no solutions for $k=2$.
It follows from a standard result in analytic number theory (as we shall see later)
that the set of $d$ for which  there is a solution with $k=2$ has natural
density $0$.
We prove the following generalization to all even exponents $k$.
\begin{thm}\label{thm:main}
Let $k\ge 2$ be even and let $r$ be a non-zero integer.
Write $\cA_{k,r}$ for the set of integers $d \ge 2$ such
that the equation
\begin{equation}\label{eqn:main}
x^k+(x+r)^k+\cdots+(x+(d-1)r)^k=y^n,\qquad x,~y,~n \in \Z, \qquad n \ge 2
\end{equation}
has a solution $(x,y,n)$. Then $\cA_{k,r}$ has natural
density $0$; by this we mean
\[
\lim_{X \rightarrow \infty} \frac{ \# \{d \in \cA_{k,r} \; :\;  d \le X\}}{X} =0.
\]
\end{thm}
If $k$ is odd, then $\cA_{k,r}$ contains all of the odd $d$:
we can take $(x,y,n)=(r(1-d)/2,0,n)$. Thus the conclusion of
the theorem does not hold for odd $k$.


\section{Some Properties of Bernoulli Numbers and Polynomials}\label{sec:Bernoulli}
In this section we summarise some classical properties of 
Bernoulli numbers and polynomials. These  are found in many
references, including \cite{Dilcher}. The Bernoulli numbers
$b_k$ are defined via the expansion 
\[
\frac{x}{e^x-1}=\sum_{k=0}^\infty b_k \frac{x^k}{k!} \, .
\]
The first few bernoulli numbers are
\[
b_0=1, \quad b_1=-1/2, \quad b_2=1/6, \quad b_3=0, \quad b_4=-1/30, \quad b_5=0, \quad b_6=1/42.
\]
It is easy to show that $b_{2k+1}=0$ for all $k \ge 1$. The $b_k$ are rational numbers,
and the Von Staudt--Clausen theorem asserts that for $k \ge 2$ even
\[
b_k+\sum_{(p-1) \mid k} \frac{1}{p}  \; \in \; \Z
\]
where the sum ranges over primes $p$ such that $(p-1) \mid k$.

The $k$-th Bernoulli polynomial can be defined by
\begin{equation}\label{eqn:berndefn}
B_k(x)=\sum_{m=0}^k \binom{k}{m} b_m x^{k-m}.
\end{equation}
Thus it is a monic polynomial with rational coefficients,
and all primes appearing in the denominators are bounded by $k+1$.
It satisfies the symmetry
\begin{equation}\label{eqn:symmetry}
B_k(1-x)=(-1)^k B_k(x),
\end{equation}
the identity
\begin{equation}\label{eqn:diff}
B_k(x+1)-B_k(x)=k x^{k-1},
\end{equation}
and the recurrence
\begin{equation}\label{eqn:derivative}
B_k^\prime(x)=k B_{k-1}(x).
\end{equation}
Whilst all the above results have been known since at least the 19th century,
we also 
make use of the following far more recent and difficult theorem 
due to Brillhart \cite{Brillhart} and Dilcher \cite{Dilcher2}.
\begin{thm}[Brillhart and Dilcher]\label{thm:Dilcher}
The Bernoulli polynomials are squarefree.
\end{thm}

\subsection*{Relation to sums of powers}
\begin{lem}\label{lem:preprevious}
Let $r$ be a non-zero integer and $k$, $d \ge 1$. Then
\[
x^k+(x+r)^k+\cdots+(x+r(d-1))=
\frac{r^k}{k+1}
\left(B_{k+1}\left(\frac{x}{r}+d\right)-B_{k+1}\left(\frac{x}{r}\right) \right).
\]
\end{lem}
This formula can be found in \cite[Section 24.4]{Dilcher},
but is easily deduced from the identity \eqref{eqn:diff}.

\begin{lem}\label{lem:previous}
Let $q\ge k+3$ be a prime. Let $a$, $r$, $d$ be integers with $d \ge 2$,
and $r \ne 0$.
Suppose $q \mid d$ and $q \nmid r$. Then
\[
a^k+(a+r)^k+\cdots+(a+r(d-1))^k \equiv r^k\cdot d \cdot B_k(a/r) \pmod{q^2}.
\]
\end{lem}
\begin{proof}
By Taylor's Theorem
\[
B_{k+1}(x+d)=B_{k+1}(x)+d \cdot B_{k+1}^\prime (x) +
\frac{d^2}{2} B_{k+1}^{(2)}(x)+\cdots+\frac{d^{k+2}}{(k+2)!} \cdot B_{k+1}^{(k+2)}(x).
\]
It follows from the assumption $q \ge k+3$
that the coefficients of $B_{k+1}$
are $q$-adic integers. Thus the coefficients of the polynomials
$B_{k+1}^{(i)}(x)/{i!}$ are also $q$-adic integers. As $q \mid d$ 
and $q \nmid r$ we have
\[
B_{k+1}\left(\frac{a}{r}+d\right)-
B_{k+1}\left(\frac{a}{r}\right) 
\equiv d \cdot B_{k+1}^\prime(a/r) \pmod{q^2}.
\]
The lemma follows from \eqref{eqn:derivative} and Lemma~\ref{lem:preprevious}.
\end{proof}

\begin{lem} \label{lem:criterion}
Let $k$, $r$ be integers with $k \ge 2$ and $r \ne 0$. 
Let $q \ge k+3$ be a prime not dividing $r$ such that the
congruence $B_k(x) \equiv 0 \pmod{q}$ has no solutions.
Let $d$ be a positive integer
such that $\ord_q(d)=1$. Then
equation \eqref{eqn:main}
has no solutions (i.e. $d \notin \cA_{k,r}$).
\end{lem}
\begin{proof}
Suppose $(x,y,n)=(a,b,n)$ be a solution to \eqref{eqn:main}.
By Lemma~\ref{lem:previous},
\[
r^k \cdot d \cdot B_k(a/r) \equiv b^n \pmod{q^2}.
\] 
However, the hypotheses of the lemma ensure that the left-hand
side has $q$-adic valution $1$. Thus $\ord_q(b^n)=1$ giving
a contradiction.
\end{proof}

\medskip

\noindent \textbf{Remarks.}
\begin{itemize}
\item For $k \ge 3$ odd, the $k$-th Bernoulli polynomial has known rational
roots $0$, $1/2$, $1$. Thus the criterion in the lemma fails to hold
for a single prime $q$. We shall in fact show that for even $k \ge 2$
there is a positive density of primes $q$ such that $B_k(x)$ has no
roots modulo $q$.
\item The second Bernoulli polynomial is $B_2(x)=x^2-x+1/6$.
This has a root modulo $q \nmid 6$ if and only if $q \equiv \pm 1 \pmod{12}$.
We thus recover the result of Bai and Zhang mentioned in 
the introduction: if $q \equiv \pm 5 \pmod{12}$ and $\ord_q(d)=1$
then \eqref{eqn:premain} has no solutions with $k=2$.
\end{itemize}

\section{A Galois property of even Bernoulli polynomials}

\begin{prop}\label{prop:Galois}
Let $k \ge 2$ be even, and let $G$ be the Galois group
of the Bernoulli polynomial $B_k$. Then there is an element
$\mu \in G$ that acts freely on the roots of $B_k$.
\end{prop}

There is a long-standing conjecture that the even Bernoulli polynomials
are irreducible; see for example \cite{Brillhart}, \cite{Carlitz}, \cite{Kimura}.
One can easily deduce Proposition~\ref{prop:Galois} from this conjecture.
We give an unconditional proof of Proposition~\ref{prop:Galois} in Section~\ref{sec:Galois}.
As noted previously, if $k$ is odd, then $B_k$ has rational roots $0$, $1/2$, $1$, so
the conclusion of the proposition certainly fails for odd $k$.

\subsection*{A Density Result}
Let $\cA$ be a set of positive integers. For $X$ positive,
define
\[
\cA(X)=\# \{d \in \cA : d \leq X \}.
\]
The \textbf{natural density} of $\cA$ is defined as the limit (if it exists)
\[
\delta(\cA)=\lim_{X \rightarrow \infty} \frac{\cA(X)}{X}.
\]
For a given prime $q$, define
\[
\cA^{(q)}=\{d \in \cA \; : \; \ord_q(d)=1\}.
\]
We shall need the following result of Niven \cite[Corollary 1]{Niven}.
\begin{thm}[Niven]\label{thm:Niven}
Let $\{q_i\}$ be a set of primes such that $\delta(\cA^{(q_i)})=0$
and $\sum q_i^{-1}=\infty$. Then $\delta(\cA)=0$.
\end{thm}

\subsection*{Proposition~\ref{prop:Galois} implies Theorem~\ref{thm:main}}
We now suppose Proposition~\ref{prop:Galois} and use it to deduce
Theorem~\ref{thm:main}.
Let $k \ge 2$ be an even integer. 
Write $G$ for the Galois group of the Bernoulli polynomial $B_k$.
Let $\mu \in G$ be the element acting freely on the roots of $B_k$
whose existence is asserted by Proposition~\ref{prop:Galois}.
By the Chebotarev density theorem  \cite[Chapter VIII]{CF} there
is a set of primes $\{q_i\}_{i=1}^\infty$ having positive Dirichlet density such that 
for each $q=q_i$, the Frobenius element $\Frob_q \in G$ is conjugate to $\mu$.
We omit from $\{q_i\}$ (without affecting the density)
the following:
\begin{itemize}
\item primes $q \le k+2$;
\item primes $q$ dividing $r$;
\item primes $q$ dividing the numerator of the discriminant of $B_k$
(which is non-zero by Theorem~\ref{thm:Dilcher}).
\end{itemize}
As $\mu$ acts freely on the roots of $B_k$,
it follows that the polynomial $B_k(x)$ has no roots modulo any of the
$q_i$.
Now let $\cA=\cA_{k,r}$ be as in the statement of Theorem~\ref{thm:main}. 
By Lemma~\ref{lem:criterion}, if $\ord_{q_i}(d)=1$ then
$d \notin \cA$. It follows that $\cA^{(q_i)}=\emptyset$.
By Theorem~\ref{thm:Niven}, we have $\delta(\cA)=0$
as required.

\section{The $2$-adic Newton polygons of even Bernoulli polynomials}
\begin{lem}\label{lem:NP}
Let $k \ge 2$ be even and write $k=2^s t$ where $t$ is odd and $s \ge 1$.
The $2$-adic Newton polygon of $B_k$ consists two segments:
\begin{enumerate}
\item[(i)] a horizantal segment joining the points $(0,-1)$ and $(k-2^s,-1)$;
\item[(ii)] a segment joining the points $(k-2^s,-1)$ and $(k,0)$ of slope $1/2^s$.
\end{enumerate}
\end{lem}
\begin{proof}
Consider the definition of $B_k$ in \eqref{eqn:berndefn}. We know that $b_0=1$, $b_1=-1/2$
and $b_m=0$ for all odd $m \ge 3$. From the Von Staudt--Clausen theorem, we know that
$\ord_2(b_m)=-1$ for even $m \ge 2$. It follows that the Newton polygon is bounded below
by the Horizontal line $y=-1$. 

We shall need to make use of the following result of Kummer 
(see \cite{Granville}): 
if $p$ is a prime, and $u$, $v$ are positive integers then
\[
\binom{u}{v}
\equiv
\binom{u_0}{v_0} \binom{u_1}{v_1} \pmod{p},
\]
where $u_0$, $u_1$ are respectively the remainder and quotient on dividing $u$ by $p$,
and likewise $v_0$, $v_1$ are respectively 
the remainder and quotient on dividing $v$ by $p$.
Here we adopt the convention $\binom{r}{s}=0$ if $r<s$. 
Applying this with $p=2$ we see that
\[
\binom{k}{2^s} =
\binom{2^s t}{2^s} 
\equiv \binom{t}{1} \equiv t \equiv 1 \pmod{2}.
\]
Thus the coefficient of $x^{k-2^s}$ in $B_k$ has $2$-adic valuation $-1$. 
Since the constant coefficient of $B_k$ also has valuation $-1$, we obtain the segment (i)
as part of the Newton polygon. We also see that for $0<v<2^s$,
\[
\binom{k}{v} \equiv 0 \pmod{2},
\]
and so the valuation of the coefficient of $x^{k-v}$ is $ \ge 0$. Finally the coefficient
of $x^k$ is $b_0=1$ and so has valuation $0$. This gives segment (ii) and completes the proof.
\end{proof}

\bigskip

\noindent \textbf{Remark}. Inkeri \cite{Inkeri} showed that $B_k$ has no rational roots
for $k$ even. His proof required very precise (and difficult) estimates for the real
roots of $B_k$. Lemma~\ref{lem:NP} allows us to give a much simpler proof
of the following stronger results. 
\begin{thm}\label{thm:twoadic}
Let $k$ be even. Then
$B_k$ has no roots in $\Q_2$.
\end{thm}
\begin{proof}
Indeed, suppose $\alpha \in \Q_2$ is a root of $B_k$. From the slopes of the Newton
polygon segments we see that $\ord_2(\alpha)=0$ or $-1/2^s$. As $\ord_2$ takes only
integer values on $\Q_2$, we see that $\ord_2(\alpha)=0$ and so $\alpha \in \Z_2$.
Let $f=2B_k \in \Z_2[x]$. Thus $f(\alpha)=0$ and so $f(\overline{\alpha})=\overline{0} \in \F_2$.
However, $\overline{\alpha} \in \F_2=\{\overline{0},\overline{1}\}$.
Now
$f(\overline{0})=\overline{(2b_k)}=\overline{1}$, and from \eqref{eqn:symmetry}
we know that $f(\overline{1})=f(\overline{0})=\overline{1}$.  This gives a contradiciton.
\end{proof}
Although Theorem~\ref{thm:twoadic} is not needed by us, its proof helps motivate
part of the proof of Proposition~\ref{prop:Galois}.

\section{Completing the proof of Theorem~\ref{thm:main}}\label{sec:Galois}
\subsection*{A little group theory}

\begin{lem}\label{lem:group}
Let $H$ be a finite group acting transitively on a finite
set $\{ \beta_1,\dotsc,\beta_n\}$. Let $H_i \subseteq H$ 
be the stabilizer of $\beta_i$, and suppose $H_1=H_2$.  
Let $\pi \; : \; H \rightarrow C$
be a surjective homomorphism from $H$ onto a cyclic group $C$.
Then there is some $\mu \in H$ acting freely on 
$\{\beta_1,\dotsc,\beta_n\}$ 
such that $\pi(\mu)$ is a generator of $C$.
\end{lem}
\begin{proof}
Let $m=\#C$ and write
$C=\langle \sigma \rangle$.
Consider the subset
\[
C^\prime=\{ \sigma^r \; : \; \gcd(r,m)=1\};
\]
this is the set of elements that are cyclic generators of $C$, and has
cardinality $\varphi(m)$, where $\varphi$
is the Euler totient function. As $\pi$ is surjective we see that
\begin{equation}\label{eqn:invCprime}
\# \pi^{-1}(C^\prime)=\frac{\varphi(m)}{m} \cdot \# H.
\end{equation}

As $H$ acts transitively on the $\beta_i$, the stabilizers $H_i$
are conjugate and so have the same image $\pi(H_i)$ in $C$.
If this image is a proper subgroup of $C$, then take $\mu$
to be any preimage of $\sigma$. Thus $\pi(\mu)=\sigma$
is a generator of $C$, and moreover, $\mu$ does not belong to
any of the stabilizers $H_i$ and so acts freely on $\{\beta_1,\dotsc,
\beta_n\}$, completing the proof in this case. Thus we
suppose that $\pi(H_i)=C$ for all $i$. It follows that
\begin{equation}\label{eqn:size}
\# \pi^{-1} (C^\prime) \cap H_i=  \frac{\varphi(m)}{m} \cdot \# H_i
= \frac{\varphi(m)}{m} \cdot \frac{\# H}{n},
\end{equation}
where the second equality follows from the Orbit-Stabilizer Theorem.
The lemma states that there is some element $\mu$ belonging to
$\pi^{-1}(C^\prime)$ but not to $\cup H_i$. Suppose otherwise.
Then 
\[
\pi^{-1}(C^\prime) \subseteq \bigcup_{i=1}^n H_i,
\]
and therefore
\begin{equation}\label{eqn:union}
\pi^{-1} (C^\prime)=\bigcup_{i=1}^n \pi^{-1}(C^\prime) \cap H_i.
\end{equation}
Now \eqref{eqn:invCprime},
\eqref{eqn:size} and \eqref{eqn:union}
together imply that the $\pi^{-1}(C^\prime) \cap H_i$ are pairwise
disjoint. This contradicts the hypothesis that $H_1=H_2$
completing the proof.
\end{proof}

\subsection*{Proof of Propostion~\ref{prop:Galois}}
We now complete the proof of Theorem~\ref{thm:main} by proving Proposition~\ref{prop:Galois}.
Fix an even $k \ge 2$, and let $L$ be the splitting field of $B_k$. Let $G=\Gal(L/\Q)=\Gal(B_k)$
be the Galois group of $B_k$. Let $\mP$ be a prime of $L$ above $2$. The $2$-adic valuation $\ord_2$
on $\Q_2$ has a unique extension to $L_\mP$ which we continue to denote by $\ord_2$. We let
$H=\Gal(L_\mP/\Q_2) \subseteq G$ be the decomposition subgroup corresponding to $\mP$.

From Lemma~\ref{lem:NP} we see that $B_k$ factors as $B_k(x)=g(x)h(x)$ over $\Q_2$ where
the factors $g$, $h$ correspond respectively to the segments (i), (ii) in the lemma.
Thus $g$, $h$ have degree $k-2^s$ and $2^s$ respectively. We denote the roots 
of $g$ by $\{\alpha_1,\dotsc,\alpha_{k-2^s}\} \subset L_\mP$ and the roots of $h$ by
$\{\beta_1,\dotsc,\beta_{2^s}\} \subset L_\mP$.  From the slopes of the segments we
see that $\ord_2(\alpha_i)=0$ and $\ord_2(\beta_j)=-1/2^s$. It clearly follows that $h$
is irreducible and therefore that $H$ acts transitively on the $\beta_j$. Moreover,
from the symmetry \eqref{eqn:symmetry} we see that $1-\beta_1$ is a root of $B_k$,
and by appropriate relabelling we can suppose that $\beta_2=1-\beta_1$. In the notation
of Lemma~\ref{lem:group}, we have $H_1=H_2$. Now let $C=\Gal(\F_\mP/\F_2)$, where $\F_\mP$
is the residue field of $\mP$. This group is cyclic generated by the Frobenius map: $\overline{\gamma}
\mapsto \overline{\gamma}^2$. 
We let $\pi \; : \; H \rightarrow C$ be the induced surjection. By Lemma~\ref{lem:group}
there is some $\mu \in H$ that acts freely on the $\beta_i$ and such that
$\pi(\mu)$ generates $C$. To complete the proof of Proposition~\ref{prop:Galois}
it is enough to show that $\mu$ also acts freely on the $\alpha_i$. Suppose
otherwise, and let $\alpha$ be one of the $\alpha_i$ that is fixed by $\mu$.
As $\ord_2(\alpha)=0$, we can write $\overline{\alpha} \in \F_\mP$
for the reduction of $\alpha$ modulo $\mP$. Now $\alpha$ is fixed by $\mu$,
and so $\overline{\alpha} \in \F_\mP$ is fixed by $\langle \pi(\mu) \rangle=C$.
Thus $\overline{\alpha} \in \F_2$ and so $\overline{\alpha}=\overline{0}$ or $\overline{1}$.
Now let $f=2 B_k(x) \in \Z_2[x]$. Thus $f(\overline{\alpha})=\overline{0}$. But
$f(\overline{0})=\overline{(2b_k)}=\overline{1}$, and from \eqref{eqn:symmetry}
we know that $f(\overline{1})=f(\overline{0})=\overline{1}$. This contradiction completes the proof. 


\end{document}